\documentclass[american]{amsart}

\usepackage[T1]{fontenc}

\usepackage[latin1]{inputenc}

\usepackage{babel}

\makeatletter

\theoremstyle{plain}    

\newtheorem{thm}{Theorem}[section]

\numberwithin{equation}{section} 

\numberwithin{figure}{section} 

\theoremstyle{plain}    

\newtheorem{lem}[thm]{Lemma} 

\theoremstyle{plain}    

\theoremstyle{plain}    

\newtheorem{claim}[thm]{Claim} 

\theoremstyle{plain}    

\theoremstyle{remark}

\newtheorem{rem}[thm]{Remark}

\theoremstyle{remark}

\newtheorem{notation-assumptions}[thm]{Notation-Assumptions}

\usepackage[all]{xy}

\xyoption{dvips}

\CompileMatrices

\makeatother

\begin{document}

\title{Vector Fields on smooth threefolds vanishing on complete intersections}

\date{\today }

\author{Thomas Eckl}



\thanks{The author gratefully acknowledges support by the
DFG-Graduiertenkolleg ``Komplexe Mannigfaltigkeiten'' at the University of
Bayreuth.}

\address{Thomas Eckl, Institut für Mathematik, Universität
Bayreuth, 95440 Bayreuth, Germany}

\email{thomas.eckl@uni-bayreuth.de}

\urladdr{http://btm8x5.mat.uni-bayreuth.de/\~{}eckl}

\maketitle




\section{Introduction}

It is a well known fact that the existence of a vector field on a K\"ahler
manifold with a special zero locus strongly influences the geometry of the
manifold.
For example, the plurigenera and some Hodge numbers vanish 
(\cite{Kob72},\cite{CaLi73}),
and if the zero locus is projective-algebraic, the variety itself will be, too
(\cite{Hwang96}). In \cite{Wahl83} J. Wahl
proved the following 
\begin{thm} \label{Wahl-Th}
Let
$X$
be a complex projective normal variety, 
$L$
an ample line bundle,
$\dim X > 1$.
If
$H^0(\mathcal{T}_X \otimes L^{-1}) \neq 0$
then:
\begin{itemize}
\item[(i)]
$L \cong \mathcal{O}(E)$,
where the effective divisor 
$E$
is a normal variety.
\item[(ii)]
$X \cong \mathrm{Proj\ }A[t]$
with
$A = \bigoplus_{n=0}^\infty H^0(E, \mathcal{O}_E(nE))$.
$X$
is the cone over 
$E$,
and~$t$
has weight 1 and
$E$
is the divisor at
$\infty$
($t=0$).
\end{itemize}
In particular every smooth complex projective variety with a vector field 
vanishing on an ample divisor is isomorphic to
$\mathbb{P}^N$.
\end{thm}
The aim of this paper is to look for similar statements in case of zero loci in
higher codimensions, with ample normal bundle for example. The main result of 
this paper is in dimension 3:
\begin{thm} \label{Q3-Theo}
Let 
$X$
be a smooth complex projective threefold. Let
$D_1, D_2$
be two ample effective divisors on
$X$
such that the scheme theoretic intersection
$C = D_1 \cap D_2$
is an irreducible reduced curve. Let
$v \in H^0(X, T_X)$
be a vector field vanishing on 
$C$.
Then
$X$
is isomorphic to
$\mathbb{P}^3$
or to the three-dimensional quadric
$Q_3$.
\end{thm}
Note that the irreducibility assumption on
$C$
is necessary: Let
$X \cong \mathbb{P}^2 \times \mathbb{P}^1$
and
$p_1:X \rightarrow \mathbb{P}^2, p_2:X \rightarrow \mathbb{P}^1$
the projections, 
$p \in L \subset \mathbb{P}^2$
a point in a line in
$\mathbb{P}^2$. 
Let
$D_1 \subset \mathbb{P}^2 \times \mathbb{P}^1$
be the blow up of
$\mathbb{P}^2$
in
$p$, 
embedded in
$X$.
There is a point 
$q \in \mathbb{P}^1$
such that the fibre
$p_2^{-1}(q)$
contains the strict transform of
$L$. 
Let 
$D_2 \subset \mathbb{P}^2 \times \mathbb{P}^1$
be the union
$p_2^{-1}(q) \cup p_1^{-1}(L)$.
Then
$D_1,D_2$
are ample divisors of
$X$
and
$C = D_1 \cap D_2$
is the union of two lines. Choose homogeneous coordinates 
$(x_0:x_1:x_2)$
on
$\mathbb{P}^2$
such that
$L = \{x_0 = 0 \}$. Then the vector field
$v = x_0\frac{\partial}{\partial x_1} \in H^0(\mathbb{P}^2, T_{\mathbb{P}^2})$
vanishes on 
$L$
and the pullback
$p_1^\ast v \in H^0(X,T_X)$
vanishes on 
$C$.

The proof of the theorem uses a modification of Wahl's theorem and Lefschetz' 
hyperplane theorem for an inductive argument to conclude that
$\mathrm{Pic}(X) = \mathbb{Z}$.
Since there exists a vector field with nonempty zero locus 
$X$
is a smooth Fano threefold. The geometry of the ample divisors discovered by 
Wahl's theorem rules out most of the possibilities in Iskhovskhi's classification
(\cite{IskI,IskII}). The remaining cases are dealt with by looking at their 
vector fields.

\subsubsection*{Notation }

Throughout the present work, let 
$X$ 
be a complex projective variety. The tangent sheaf
$\mathcal{T}_X = \mathcal{H}om_{\mathcal{O}_X}(\Omega^1_X,\mathcal{O}_X)$
is the dual of the sheaf of differentials
$\Omega^1_X$
on
$X$. 
If
$X$
is smooth,
$T_X$
is locally free. A vector field
$v \in H^0(X,\mathcal{T}_X)$
is a global section of 
$\mathcal{T}_X$.
The groups
$\mathbb{C}^\ast, \mathbb{C}^+$
denote the multiplicative and the additive group of complex numbers.

\subsubsection*{Acknowledgement }

This paper is part of a PhD thesis written while the author was a member of the
DFG-Graduiertenkolleg ``Komplexe Mannigfaltigkeiten'' at the University of 
Bayreuth. The author is grateful to his advisor Th. Peternell for pointing out
the problem and continuous encouragement to work on it. He would like to 
thank the other members of the Graduiertenkolleg, especially S.~Kebekus and 
H.~Chr.~von~Bothmer, for creating a stimulating athmosphere and many helpful 
discussions on the topic.

\section{A Theorem of Wahl}

An easy corollary of Theorem~\ref{Wahl-Th} is the following result already 
proven in \cite{MS78}:
\begin{thm}
Let
$X$
be a complex projective smooth variety,
$L$
an ample line bundle and
$T_X$
the tangent bundle. If there is a vector field 
$v \in H^0(X,T_X)$
vanishing on an ample effective divisor 
$D \in H^0(X,L)$
then
\[ (X,L) = (\mathbb{P}^n,\mathcal{O}(1))\ \mathrm{or}\ 
           (\mathbb{P}^1,\mathcal{O}(2)). \]
\end{thm} 
The proof of Theorem~\ref{Wahl-Th} uses the normality of
$X$
when applying the Kodaira vanishing theorem. A thorough analysis of this proof 
and the ideas in \cite{MS78} shows that one can replace normality by the 
existence of an effective divisor in the ample line bundle to get a slightly
weaker result:
\begin{thm} \label{var-Wahl-Th}
Let
$X$
be a complex projective variety, 
$\dim X > 1$,
and
$L$
an ample line bundle with
$H^0(X,L) \neq 0$.
If
$H^0(X,\mathcal{T}_X \otimes L^{-1}) \neq 0$
then there is a finitely generated graded
$\mathbb{C}$-algebra
$A \subset \bigoplus_{n=0}^\infty H^0(X,nL)$
and a homogeneous element
$T \in H^0(X,L)$
such that
\[ A[T] = \bigoplus_{n=0}^\infty H^0(X,nL). \]
In this case
$X$
is isomorphic to a cone over an ample Divisor
$E$
of
$\mathrm{Proj}\ A$.
\end{thm}

\begin{proof}
The last part is an easy consequence of the first statements, 
cf.~\cite[5.3]{BS95}.

The first part can be proven exactly as in \cite{Wahl83} provided that one is
still able to construct a derivation of weight -1 on
$R = \bigoplus_{n=0}^\infty H^0(X,nL)$:
By assumption there is an effective divisor
$D \in H^0(X,L)$
and a vector field
$v \in H^0(X,\mathcal{T}_X)$
such that
$v_{|D} = 0$.
Let
$G \subset \mathrm{Aut\ }X$
be the subgroup of automorphisms fixing
$D$. 
The existence of 
$v$
implies that
$G$
is nontrivial. Furthermore the linear representation of
$G$
on the 
$m$-jets 
of a point
$x \in D$
is faithful if
$m \gg 0$. 
Therefore 
$G$
is linear algebraic and contains a linear algebraic one parameter group
$H$,
i.e.
$H \cong \mathbb{C}^\ast$
or
$\mathbb{C}^+$.

Since
$H$
stabilizes the divisor
$D$
there is a
$H$-linearization 
of the line bundle
$L$
and a dual 
$H$-action 
$\tau$
on
$R = \bigoplus_{n=0}^\infty H^0(X,L^n)$
(cf. \cite[Prop. 1.5]{GIT}). If
$H = \mathbb{C}^\ast$
there are semi-invariant elements
$F_0,F_1,\ldots,F_m$
of
$\tau$, 
i.e.
\[ \tau(t) \cdot F_i = t^{\chi(F_i)} \cdot F_i \]
which generate
$R$
as a
$\mathbb{C}$-algebra (
$\chi(F)$ 
is the weight of the semi-invariant element 
$F$).
Let
$F_0 \in H^0(X,L)$
correspond to 
$D$.

Let 
$R/F_0 = \bigoplus_m H^0(D,mL_{|D})$
be the homogeneous coordinate ring of
$D$.
Since
$\mathbb{C}^\ast$
acts trivial on
$D$
via
$\tau$
it acts semi-trivial on
$H^0(D,mL_{|D})$
with weight
$\chi_m$,
and the quotients
$\frac{\chi_m}{m}$
are equal for all
$m \geq 1$.
This implies
\[ s := \frac{\chi(F_1)}{\deg F_1} = \ldots = \frac{\chi(F_m)}{\deg F_m}, \]
and 
$ s \neq \frac{\chi(F_0)}{\deg F_0}$
since otherwise
$\mathbb{C}^\ast$
would act trivial on all of
$X$.

Now one can twist the
$\mathbb{C}^\ast$-linearization: 
Let
$\mathbb{C}^\ast$
act on
$R$
via
\[ \sigma(t) \cdot F = t^{- s\deg F} \cdot \tau(t) F \]
for all homogeneous
$F \in R$.
Hence for an arbitrary polynomial
$P \in \mathbb{C}[X_0,\ldots,X_m]$:
\[ \sigma(t) P(F_0,F_1,\ldots,F_m) = P(t^{\chi(F_0)-s}F_0,F_1,\ldots,F_m), \]
and
$\chi(F_0)-s \neq 0$.
It follows 
$D_\sigma F_1 = \ldots = D_\sigma F_m = 0, D_\sigma F_0 = (\chi(F_0)-s)F_0$
where 
$D_\sigma$
is the R-derivation corresponding to
$\sigma$.
One can divide 
$D_\sigma$
by
$F_0$
to get a
$(-1)-\!$
derivation
$D_{-1}$.

If
$H = \mathbb{C}^+$
let
$F_0, F_1, \ldots, F_m$
still be homogeneous generators of 
$R$, 
the element
$F_0$
corresponding to 
$D$.
Since the unipotent group
$\mathbb{C}^+$
fixes 
$D = \mathrm{Proj}(R/F_0)$
the group acts trivial on
$R/F_0$.
This implies for homogeneous
$F$
that
\[ \tau(t) F - F \in (F_0). \]
Once more one can divide the corresponding derivation
$D_\tau$
by
$F_0$
and gets a
$(-1)-\!$
derivation
$D_{-1}$
on
$R$.

Now one constructs an element 
$t \in H^0(X,L)$
as in \cite[Lemma 2.7]{Wahl83} with
\[ D_{-1}t = 1. \]
Then
$R \cong A[T]$
with
$A = \{ r \in R | D_{-1}r = 0 \}$
and
$\deg T = 1$
(cf. \cite[Prop. 2.4]{Wahl83}),
and the theorem follows.
\end{proof}

\section{Reduction to Fano Manifolds with Picard number 1} \label{Red-Sec}

As in Theorem~\ref{Q3-Theo} let 
$D_1, D_2$
be two ample effective divisors on a smooth complex projective variety
$X$
of dimension 3 such that the scheme theoretic intersection 
$C = D_1 \cap D_2$
is an irreducible reduced curve. This implies that
$D_1,D_2$
are irreducible and reduced.

Let furthermore 
$v \in  H^0(X,T_X)$
be a vector field vanishing on
$C$.
\begin{lem} \label{Stab-Lem}
If
$D_2$
is not stabilized by 
$v$
then there will be an irreducible and reduced divisor
$D_\infty$
stabilized by
$v$,
which is linearly equivalent to
$D_1$
and whose scheme theoretic intersection with 
$D_2$
equals
\[ D_\infty \cap D_2 = D_1 \cap D_2 = C. \]
\end{lem}
\begin{proof}
Let
$G \subset \mathrm{Aut\ }X$
be the connected and nontrivial algebraic subgroup of the automorphism group of
$X$
which fixes the zero locus 
$Z(v)$
of
$v$. 
Because the action of 
$G$
on the vector space of
$m$-jets 
at a fixed point will be faithful for
$m \gg 0$,
the group 
$G$
is linear algebraic. 

Let
$H \subset G$
be the minimal algebraic subgroup whose Lie algebra contains
$v$.
A representation of 
$G$
in
$GL(V)$
shows that
$H \cong (\mathbb{C}^\ast)^k \times (\mathbb{C}^+)^l$
(\cite[II.7.3]{BorelLAG}). Because 
$H$
is commutative, the fixed point locus
$X^H$
is contained in
$Z(v)$
and
$H$
stabilizes
$Z(v)$.
By composing the various 
$\mathbb{C}^\ast$-
and
$\mathbb{C}^+$-actions 
one can move the divisor 
$D_1$
along orbits to a linearly equivalent divisor
$D_\infty$
stabilized by 
$H$.

$H(Z(v)) = Z(v)$
implies
$D_1 \cap D_2 \subset D_\infty \cap D_2$,
linear equivalence means 
$D_1.D_2 = D_\infty.D_2$,
consequently
$D_1 \cap D_2 = D_\infty \cap D_2$.
And
$D_\infty$
is reduced and irreducible because 
$D_\infty \cap D_2$
is. 
\end{proof}

Assume from now on that 
$D_1$
is stabilized by the vector field
$v$.
If
$v_{|D_1} = 0$,
Wahl's Theorem~\ref{Wahl-Th} will imply
$X \cong \mathbb{P}^3$.
If 
$v_{|D_1} \neq 0$,
the variant Theorem~\ref{var-Wahl-Th} will imply that
$D_1$
is a cone
$\mathcal{C}(C,L)$
for a (possibly singular) curve
$C$
and an ample line bundle 
$L$
on 
$C$.

The cone
$D_1 \cong \mathcal{C}(C,L)$
is the contraction of the section of the projective space bundle
$\mathbb{P}(\mathcal{O}_C \oplus \mathcal{O}_C(L))$
belonging to the projection
$\mathcal{O}_C \oplus \mathcal{O}_C(L) \rightarrow \mathcal{O}_C$.
Therefore, 
$H^2(D_1, \mathbb{Q}) = \mathbb{Q}$,
and the finitely generated abelian group 
$H^2(D_1, \mathbb{Z})$
has rank 1. Now apply Lefschetz' hyperplane theorem: the natural map
\[ H^2(X, \mathbb{Z}) \hookrightarrow H^2(D_1, \mathbb{Z}) \]
is an injection. Since 
$\mathrm{Num}(X)$
is a torsion free quotient of 
$\mathrm{NS}(X) \subset H^2(X, \mathbb{Z})$,
it follows that
$\mathrm{Num}(X) = \mathbb{Z}$. 
There is an ample divisor
$H$
on
$X$
and
$r \in \mathbb{Z}$
such that
\[ K_X \equiv r H\]

On the other hand
$X$
can be covered by rational curves: Since there is a non trivial vector field 
with zeroes on
$X$
one of the groups
$\mathbb{C}^\ast$
or
$\mathbb{C}^+$
is acting on
$X$
(s. proof of Lemma~\ref{Stab-Lem}). The closures of the orbits are rational 
curves. Consequently,
$K_X$
is not nef (\cite[II.3.13.1]{RC}), not ample, and
$-K_X$
is ample. By \cite[Prop.1.15]{IskI} (and \cite{Shok80} for the proof of 
hypothesis~1.14 in \cite{IskI}),
$X$
is a smooth Fano threefold with
\[ \mathrm{Pic} X \cong H^2(X, \mathbb{Z}) = \mathbb{Z}. \]
A classification of these Fano threefolds is given by the following table
(cf.~\cite{IskI},\cite{IskII}):

\begin{center}
 \begin{tabular}{|c|r|r|r|l|} \hline
   $r$ & $(H)^3$ & $b_3/2$ & $g$ & $X$ \\ \hline
   $4$ & $1$ & $0$ & $33$ & $\mathbb{P}^3$ \\ \hline
   $3$ & $2$ & $0$ & $28$ & $Q \subset \mathbb{P}^4$, the quadric \\ \hline
   $2$ & $1$ & $21$ & $5$ & $V_1$, a covering of the cone over the
                                  Veronese surface \\
   $2$ & $2$ & $21$ & $9$ & $V_2$, a double covering of 
                                   $\mathbb{P}^3$ \\
   $2$ & $3$ & $5$ & $13$ & $V_3 \subset \mathbb{P}^4$ , a cubic \\
   $2$ & $4$ & $2$ & $17$ & $V_4 \subset \mathbb{P}^5$ , an intersection of 
                                                         two quadrics \\
   $2$ & $5$ & $0$ & $21$ & $V_5$, the intersection 
                                   $Gr(1,4) \subset \mathbb{P}^9$
                           with $\mathbb{P}^6$ \\ \hline
   $1$ & $2$ & $52$ & $2$ & $V_2^\prime$ , a double covering of
                                 $\mathbb{P}^3$ \\
   $1$ & $4$ & $30$ & $3$ & $V_4^\prime \subset \mathbb{P}^4$ , a quartic \\
   $1$ & $4$ & $30$ & $3$ & $V_4^{\prime\prime}$, a double covering
                                                  of a quadric \\ \hline
   \end{tabular}
\end{center}

\begin{center}
   \begin{tabular}{|c|r|r|r|l|} \hline
   $r$ & $(H)^3$ & $b_3/2$ & $g$ & $X$ \\ \hline
   $1$ & $6$ & $20$ & $4$ & $V_6 \subset \mathbb{P}^5$ , an intersection of a 
                               quadric with a cubic \\
   $1$ & $8$ & $14$ & $5$ & $V_8 \subset \mathbb{P}^6$, an intersection of 
                                                        three quadrics \\
   $1$ & $10$ & $10$ & $6$ & $V_{10} \subset \mathbb{P}^7$ \\
   $1$ & $12$ & $5$ & $7$ & $V_{12} \subset \mathbb{P}^8$ \\
   $1$ & $14$ & $5$ & $8$ & $V_{14} \subset \mathbb{P}^9$ \\
   $1$ & $16$ & $3$ & $9$ & $V_{16} \subset \mathbb{P}^{10}$ \\
   $1$ & $18$ & $2$ & $10$ & $V_{18} \subset \mathbb{P}^{11}$ \\
   $1$ & $22$ & $0$ & $12$ & $V_{22} \subset \mathbb{P}^{13}$ \\ \hline
   \end{tabular}
\end{center}  

In this table,
$H$
is the ample generator of
$\mathrm{Pic}(X)$,
$r$
is the index of
$V$,
i.e.
$-K_X = r \cdot H$,
and
$g = -(K_X)^3/2 + 1$
is the genus of 
$X$.

The ample divisor 
$H$
is very ample except when 
$X$
is of type 
$V_1, V_2, V_2^\prime, V_4^{\prime\prime}$.
These cases are dealt with in the last section, while in the next section
one assumes that
$H$
is very ample.

\section{An estimate for the degree} \label{Deg-Sec}

By assumption 
$D = d \cdot H$ 
and
$D_2 = d_2 \cdot H$
are very ample divisors on
$X$.
Theorem~\ref{var-Wahl-Th} implies that
$D$
is a cone
$\mathcal{C}(C,\mathcal{O}(D_{2|C}))$
over a (possibly singular) curve 
$C$
with vertex 
$P$,
defined by the divisor
$D_2$
restricted to
$C$.
\begin{claim} \label{d_2-lem}
$d_2 = 1$.
\end{claim}
\begin{proof}
By construction
$C$ 
is linearly equivalent to
$D_{2|D}$,
hence very ample. The corresponding embedding
$D \hookrightarrow X \stackrel{|D_2|}{\rightarrow} \mathbb{P}$
maps
$C$
into a hyperplane, and the cone
$D$
consists of the lines through the cone vertex
$P$
(not in the hyperplane) and points
$Q \in C$. 
Such a line 
$L$
cuts
$C$
transversally in one point, consequently (in
$X$):
\[ 1 = C . L = d_2H.L. \]
$H.L > 0$
implies the hypothesis.
\end{proof}

Let
$\pi: \widetilde{X} \rightarrow X$
be the blow-up of 
$X$
in the cone vertex
$P$
and
$E$
the exceptional divisor. Because of the universal property of the blow-up the 
strict transform
$S$
of
$D$
also is the blow-up of
$D$
in 
$P$.
By \cite[6.7.1]{FIS} the effective Cartier divisor
$S$
is linearly equivalent to
$\pi^\ast(D) - \mu_D E$
for a
$\mu_D \in \mathbb{Z}$,
the multiplicity of the point
$P$
in the variety
$D$.
\begin{claim}
$\mu_D = dH^3$.
\end{claim}
\begin{proof}
Since
$H$
is very ample, by Bertini there exist two smooth hyperplane sections 
$H_1, H_2 \in |H|$
with
$H_1 \cap H_2 \cap D = \{P\}$.
Then 
$H_1 \cap H_2$
does not contain any line
$L$
from 
$P$
to a point
$Q \in C$.

As in the proof of the previous claim the 
$H_i$
intersect every line 
$L$
from 
$P$
to a point
$Q \in C - H_i$
exactly in
$P$
with intersection multiplicity
$1$,
i.e.\/ transversal. Therefore the strikt transforms 
$\tilde{D},\tilde{H_1},\tilde{H_2}$
do not intersect at all. Since the
$H_i$
are smooth, the intersection multiplicities of the
$H_i$
in
$P$
are
$\mu_{H_i} = 1$.
Then, by \cite[12.4.8]{FIS}
\[ \mu_D = \mu_D \mu_{H_1} \mu_{H_2} = D.H_1.H_2 = dH^3. \] 
\end{proof}

The strict transform
$S$
of the cone
$D$
is isomorphic to the 
$\mathbb{P}^1-\!$
bundle
\[ S \cong \mathbb{P}(\mathcal{O}_C \oplus \mathcal{O}_C(C)) \]
over the curve
$C$.
Let
$f: \widehat{C} \rightarrow C$
the normalization of the possibly singular curve
$C$ 
and let the 
$\mathbb{P}^1-\!$
bundle
\[ \widehat{S} \cong \mathbb{P}(\mathcal{O}_{\widehat{C}} \oplus 
                            f^\ast \mathcal{O}_C(C)) \]
over 
$\widehat{C}$
be the normalization of
$S$.

Thus one has the following diagram:
\[ \xymatrix{
             & {\tilde{X}} \ar[r]^{\pi} & {X} \\
            {\hat{S}} \ar[r]^{f} \ar[d]_{\hat{p}} & {S} \ar@{^{(}->}[u] 
           \ar[r]^{\pi} \ar[d]^{p} & {D} \ar@{^{(}->}[u] \\
            {\hat{C}} \ar[r]^{f} & {C}  \\
            } \]

\noindent
\begin{claim}
$H^3 \leq 4$.
\end{claim}
\begin{proof}
By the adjunction formulas one gets
\[ -K_{X|D} = -K_D + D_{|D} \]
and
\[ \begin{array}{ccl}
   K_S & = & K_{\tilde{X}|S} + S_{|S} \\ 
       & = &\pi^\ast(K_X)_{|S}+ 2 E_{|S} +\pi^\ast(D)_{|S}-\mu E_{|S} = \\
       & =  & \pi^\ast(K_{X|D} + D_{|D}) + (2 - \mu)E_{|S} =
                       \pi^\ast(K_D) + (2 - \mu)E_{|S}.
   \end{array}\]
The nonnormal locus 
$N$
on
$S$
is given by the conductor ideal of the normalization 
$f: \widehat{S} \rightarrow S$,
and its support consists of whole fibers of the 
$\mathbb{P}^1-\!$
bundle
$S \rightarrow C$.
From the subadjunction formula for normalization (and the formula for the 
canonical bundle on smooth ruled surfaces), it follows
\[ f^\ast K_S =  K_{\widehat{S}} - N  = -2 C_0 + kF,\ k \in \mathbb{Z},\]
where
$C_0 \cong \widehat{C}$
is the section with negative self intersection on the smooth ruled surface
$\widehat{S}$
and
$F$
is a fiber.

Let
$F_{\widehat{S}}$
be a fiber of
$\widehat{p}: \widehat{S} \rightarrow \widehat{C}$
and
$F_S$
a general fiber of
$p: S \rightarrow C$.
Then
$f_\ast F_{\widehat{S}} = F_S$. 
Furthermore,
$\pi_\ast F_S = L$,
where 
$L$
is a line from the cone vertex 
$P$
to a point
$Q \in C$.
Since
$E_{|S}$
is a section of the projective line bundle
$S$
over
$C$, 
the intersection multiplicity 
$F_S.E_{|S} = 1$. 
Therefore:
\[ \begin{array}{ccl}
   -2 & = & (-2C_0 + kF).F = f^\ast K_S . F_{\widehat{S}} = 
             K_S.F_S = \\
      & = & \pi^\ast(K_D).F_S + (2 - \mu)E_{|S} . F_S = 
            K_D.\pi_\ast F_S + 2- H^3 d = \\
      & = & K_D.L + 2 - H^3 d. 
   \end{array} \]
Let
$r_D$
be defined by
$-K_D = - r_D H_{|D}$.
As in the proof of Claim~\ref{d_2-lem},
$H_{|D}.L = 1$.
Consequently:
\[ -2 = - r_D H_{|D}.L + 2 - H^3 d = -r_D + 2 - H^3 d, \]
and this implies 
$r_D = 4 - H^3 d$.
But the index 
$r_X$
of the Fano threefold 
$X$
is
$\geq 1$, 
hence
\[ r_X = r_D + d = 4 - (H^3 - 1)d \geq 1, \]
and
$H^3 \leq 4$.
\end{proof}

Now, the following cases must be considered:
\begin{itemize}
\item
$H^3 = 4$:
Then
$d = 1$
and
$r_X = 1$.
Under the assumption that 
$H$
is very ample,
$X$
must be a quartic.
\item
$H^3 = 3$:
Then still
$d = 1$,
but
$r_X = 2$.
So
$X$
must be a cubic.
\item
$H^3 = 2$:
Then
$d = 1, 2$
or
$3$, corresponding to
$r_X = 3, 2, 1$,
and
$X$
must be a quadric (if
$H$
is very ample).
\item
$H^3 = 1$:
This implies 
$r_X = 4$
and
$X = \mathbb{P}^3$.
\end{itemize}

\section{Special Fanos}

In this section, vector fields vanishing on the reduced and irreducble 
intersection of two (very) ample divisors on
$\mathbb{P}^3$
and the quadric
$Q_3$
will be constructed, and it will be 
shown that such vector fields do not exist on a cubic, on a quartic and on 
varieties of type
$V_1, V_2, V_2^\prime, V_4^{\prime\prime}$,
where
$H$
is not very ample.

\subsection{$\boldsymbol{V_1, V_2, V_2^\prime}$ und  
            $\boldsymbol{V_4^{\prime\prime}}$}

These Fano varieties are described in \cite{IskII}:
\begin{itemize}
\item[(a)]
The morphism
$\varphi_{K_{V_1}^{-1}}: V_1 \rightarrow W_4 \subset \mathbb{P}^6$
induced by the complete linear system
$|K_{V_1}^{-1}|$
is a 2:1-covering of the cone
$W_4$
over the Veronese surface
$F_4 \subset \mathbb{P}^5$.
It branches over the smooth divisor
$S \subset W_4$
cut out by a cubic hypersurface not containing the cone vertex.
\item[(b)]
The morphism
$\varphi_{H}: V_2 \rightarrow \mathbb{P}^3$
is a 2:1-covering with smooth branching divisor
$S \subset \mathbb{P}^3$
of degree
$4$.
\item[(c)]
The morphism
$\varphi_{H}: V_2^\prime \rightarrow \mathbb{P}^3$
is a 2:1-covering with smooth branching divisor
$S \subset \mathbb{P}^3$
of degree
$6$.
\item[(d)]
The morphism
$\varphi_{H}: V_4^{\prime\prime} \rightarrow Q_3 \subset \mathbb{P}^4$
is a 2:1-covering of the threedimensional quadric with smooth branching divisor
$S \subset Q_3$
cut out by a hypersurface of degree 4.
\end{itemize}

Assume that on one of these Fano threefolds there exists a vector field
$v \in H^0(X,T_X)$
vanishing on a curve. Then there is a linear algebraic group
$G$
acting on
$X$
and inducing
$v$.
Since the coverings are induced by complete linear systems, 
$G$
acts on the bases 
$W$
of the coverings, too. Furthermore
$G$
stabilizes the branching divisors
$S \in W$
since these divisors describe the locus of the points where the rank of the 
differentials
$df_x: T_{X,x} \rightarrow T_{W,f(x)}$
drops.

Therefore,
$v$ 
restricted to 
$S$ 
is a (non trivial) vector field on the branching divisor. By assumption,
$v$
should vanish on a curve 
$C$. 
Since the branching divisor 
$S$
is ample,
$S$
intersects
$C$,
and
$v_{|S}$
has a zero. As explained in section~\ref{Red-Sec} this implies the existence
of a
$\mathbb{C}^+$-
or
$\mathbb{C}^\ast$-action 
on 
$S$. 
Consequently, 
$S$
is covered by rational curves, i.e. the closures of the orbits.
$S$
is uniruled, and
$K_S$
is not nef.

On the other hand it is easy to compute 
$K_S$
since the branching divisors are complete intersections:
\begin{itemize}
\item[(a)]
On varieties of type 
$V_1$
is
$K_S = \mathcal{O}_S(2)$,
because the Veronese surface and the cone over it are generated by three 
quadrics.
\item[(b)]
If
$X$
is of type
$V_2$,
then
$K_S = \mathcal{O}_S$.
\item[(c)]
$X$
of type
$V_2^\prime$:
then
$K_S = \mathcal{O}_S(2)$.
\item[(d)]
$X$
of type
$V_4^{\prime\prime}$:
then 
$K_S = \mathcal{O}_S(1)$.
\end{itemize}
In all cases
$K_S$
is nef, contradiction.

\subsection{Quadrics, cubics and quartics.}
The starting point is the following 
\begin{lem}
All vector fields on a hypersurface
$H \in \mathbb{P}^n$,
$n > 2$,
are induced by equivariant vector fields on
$\mathbb{P}^n$.
\end{lem}
\begin{proof}
Set
$d := \deg H$.
First consider the structure sheaf sequence of
$H$
tensorized by the tangent bundle
$T_{\mathbb{P}^n}$,
\[ 0 \rightarrow T_{\mathbb{P}^n}(-d) \rightarrow T_{\mathbb{P}^n} 
     \rightarrow T_{\mathbb{P}^n} \otimes \mathcal{O}_H \rightarrow 0, \]
and the beginning of the corresponding long exact sequence,
\[ 0 \rightarrow H^0(\mathbb{P}^n, T_{\mathbb{P}^n}(-d)) \rightarrow
                 H^0(\mathbb{P}^n, T_{\mathbb{P}^n}) 
                 \stackrel{\nu}{\rightarrow}
                 H^0(H, T_{\mathbb{P}^n} \otimes \mathcal{O}_H) \rightarrow
                 H^1(\mathbb{P}^n, T_{\mathbb{P}^n}(-d)). \]
The Euler sequence
\[ 0 \rightarrow \mathcal{O}_{\mathbb{P}^n}(-d) \rightarrow
   \mathcal{O}_{\mathbb{P}^n}(1-d)^{n+1} \rightarrow T_{\mathbb{P}^n}(-d)
   \rightarrow 0 \]
and
$H^i(\mathbb{P}^n, \mathcal{O}_{\mathbb{P}^n}(k)) = 0$
for
$0 < i < n$,
$k \in \mathbb{Z}$
imply
$H^1(\mathbb{P}^n, T_{\mathbb{P}^n}(-d)) = 0$.
Therefore,
$\nu$
is surjective. By the normal sequence
\[ 0 \rightarrow \mathcal{T}_H \rightarrow T_{\mathbb{P}^n} \otimes 
     \mathcal{O}_H
     \rightarrow N_{H/\mathbb{P}^n} \rightarrow 0 \]
one finally has 
$H^0(H, \mathcal{T}_H) \subset H^0(H, T_{\mathbb{P}^n} \otimes \mathcal{O}_H)$.
\end{proof}

From now on, let
$H \subset \mathbb{P}^4$
be a smooth hypersurface of degree
$d = 2,3,4$
with a vector field
$v \in H^0(H,T_H)$
vanishing on the irreducible and reduced intersection
$C = H_1 \cap H_2$
of two very ample divisors. The lemma above shows that 
$v$
is induced by a vector field on
$\mathbb{P}^4$
also called
$v$.

Claim~\ref{d_2-lem} and the cases at the end of section~\ref{Deg-Sec}
show that 
$v$
stabilizes a cone
$D = H \cap H_1$
cut out by a hyperplane, and the zero locus of
$v$
is the intersection
$C = H \cap H_1 \cap H_2$
with another hyperplane. The cone vertex
$P \in D$
is not contained in
$H_2$
since otherwise, 
$D \cap H_2$
contains more than one line because of
$\deg H > 1$. 

Now choose homogeneous coordinates
$(x_0: x_1 : x_2 : x_3 : x_4)$
on
$\mathbb{P}^4$
such that the cone
$D$
lies on the hyperplane
$H_1 = \{ x_4 = 0 \}$,
and the cone basis lies
in the plane 
$H_1 \cup H_2 = \{ x_3 = x_4 = 0 \}$.
Furthermore, let
$P =(0:0:0:1:0)$.
\begin{lem} \label{HypPol-Lem}
The hypersurface 
$H$
is given by a polynomial
\[ h = f(x_0,x_1,x_2) + x_4 g(x_0,x_1,x_2,x_3,x_4), \]
where
$\deg f = d$,
$\deg g = d-1$. 
The coefficients of the monomial
$x_3^{d-1}$
in
$g$
and of the monomial
$x_i x_4^{d-1}$
in
$h$ 
do not vanish
(at least for one 
$0 \leq i \leq 4$).
\end{lem}
\begin{proof}
The polynomial
$h$
may be written as
\[ h = f(x_0,x_1,x_2) + x_3 k(x_0,x_1,x_2,x_3) + x_4 g(x_0,x_1,x_2,x_3,x_4). \]
Since the intersection with
$H_1 = \{ x_4 = 0\}$
is supposed to be a cone over a basis
$C \subset \{ x_3 = x_4 = 0 \}$
with vertex
$(0:0:0:1:0)$,
this cone
$D$
is given by the equation
$f(x_0,x_1,x_2) = 0$,
and
$k \equiv 0$.

Since the hypersurface is smooth, the rest follows from computing the 
gradient of
$h$:
the coefficient of
$x_3^{d-1}$
must not be
$0$
because otherwise the vertex
$(0:0:0:1:0)$
will be a singularity in
$H$,
too.
Similarly, the coefficient of at least one of the monomials
$x_i x_4^{d-1}$,
$i \leq 4$,
must not vanish,
because otherwise the point
$(0:0:0:0:1)$
lies in
$H$
and will be not smooth.
\end{proof}

Now, by the Euler sequence vector fields on
$\mathbb{P}^4$
correspond to homogeneous derivations
$D = \sum_{i,j} a_{ij} x_i \frac{\partial}{\partial x_j}$
of weight
$0$
on
$\mathbb{C}[x_0, \ldots , x_4]$,
modulo the Euler vector fields
$\sum_i a x_i \frac{\partial}{\partial x_i}$.
\begin{lem}
The vector field 
$v$
corresponds to the derivation 
$D_v = \sum_{i,j} a_{ij} x_i \frac{\partial}{\partial x_j}$
given by (a non trivial scalar multiple of) the matrix
\[ \left( a_{ij} \right) =
   \left( \begin{array}{cccc|c}
      0 & \cdots & \cdots & 0 & 0 \\ 
      \vdots & 0 & & \vdots & \vdots \\
      \vdots & & 0 & 0 & \vdots \\
      0 & \cdots & 0 & 1 & 0 \\ \hline
      \ast & \ast & \ast & \ast & a
   \end{array} \right). \]
\end{lem}
\begin{proof}
Since $\deg C > 1$,
the curve
$C$
is no line.
Since vector fields on projective spaces always vanish on complete linear 
subspaces,
$v$
vanishes on the plane
$x_3 = x_4 = 0$.
On the other hand, 
$v$
does not vanish on the hyperplane
$x_4 = 0$,
because it contains the cone
$D$,
and
$v$
should be non trivial on
$D$.
Since every vector field on a cone vanishes in the vertex,
$v$
vanishes in
$P = (0:0:0:1:0)$.

This implies the hypothesis, because the zero locus of vector fields on
$\mathbb{P}^4$
consists of the
eigenspaces of the transposed corresponding matrix.
\end{proof}

Now,
$v$
stabilizes the hypersurface
$H$
exactly when the derivation
$D_v$
maps the principal ideal
$(h) \subset \mathbb{C}[x_0, \ldots , x_4]$
describing 
$H$
to itself, i.e.
$D_v h = \lambda h$
for a
$\lambda \in \mathbb{C}$.
Since
$a_{ij} = 0$
for
$i,j = 0,1,2$,
\[ D_v h = D_v f + D_v x_4 \cdot g + x_4 D_v g = 0 + x_4 (ag + D_v g), \]
and consequently
$\lambda = 0$.

This gives immediately a quadric with a vector field vanishing as in 
Theorem~\ref{Q3-Theo}: Let 
$h = x_0^2 + x_1^2 + x_2^2 + x_3 x_4$
be the equation of the matrix and let
$v$ 
correspond to the derivation
$D_v = x_3 \frac{\partial}{\partial x_3} - x_4 \frac{\partial}{\partial x_4}$.
Then,
$D_v h = 0$,
and
$v$
vanishes on the smooth quadric
$x_0^2 + x_1^2 + x_2^2 = 0$
contained in the plane
$x_3 = x_4 = 0$
and in the cone vertex
$(0:0:0:1:0)$.
Furthermore,
$v$
stabilizes the cone
$x_0^2 + x_1^2 + x_2^2 = 0$
contained in the hyperplane
$x_4 = 0$.

Why are there no such vector fields on cubics and quartics ?
\begin{itemize}
\item[(a)]
By Lemma~\ref{HypPol-Lem} the monomial
$x_3^{d-1} x_4$
has a coefficient
$c \neq 0$
in 
$h$,
but there is no monomial of the form
$x_3^{d-1} x_i$,
$i \leq 3$,
in 
$h$.
\item[(b)]
The coefficient of
$x_3^{d-1} x_4$
in
$D_v h = 0$
is
$0$.
Decomposing
$h$
in monomials and
$D_v$
in ``monomial'' derivations of the form
$a_{ij}x_i\frac{\partial}{\partial x_j}$,
one sees that only the derivations
\[ a_{3j}x_3\frac{\partial}{\partial x_j} x_jx_3^{d-2}x_4 ,\
   a_{4j}x_4\frac{\partial}{\partial x_j} x_jx_3^{d-1}x_j \]
contribute to the coefficient of
$x_3^{d-1} x_4$.
Since
$a_{3j} = 0$
for 
$i = 0,1,2,4$
and the coefficients of
$x_3^{d-1} x_i$,
$i \leq 3$,
vanish, the coefficient of
$x_3^{d-1} x_4$
in 
$D_v h$ 
is the coefficient of
\[ c(a_{33}x_3\frac{\partial}{\partial x_3} x_3^{d-1}x_4 +
     a_{44}x_4\frac{\partial}{\partial x_4} x_3^{d-1}x_4), \]
i.e.
$c(d-1 +a)$.
\item[(c)]
Consequently, 
$a = 1 - d$. 
This number is an eigenvalue of the matrix corresponding to the derivation,
different from the other Eigen values.
Therefore, in appropriate coordinates the matrix is diagonal, i.e.
\[ D_v = x_3\frac{\partial}{\partial x_3} + 
        (1-d)x_4\frac{\partial}{\partial x_4}. \] 
\item[(d)]
The difference between quadrics and cubics resp. quartics comes from the fact 
that the monomials 
$x_3 x_4^{d-1}$
and 
$x_3^{d-1} x_4$
are different for
$d > 2$. 
The coefficient
$c$
of
$x_3 x_4^{d-1}$ 
in
$h$
must be
$0$
for 
$d > 2$,
because otherwise the coefficient of
$x_3 x_4^{d-1}$
in
$Dh$
equals
$c(1-(d-1)^2) \neq 0$.
The same is true for the monomials
$x_i x_4^{d-1}$
with
$i = 0,1,2,4$.
But this a contradiction to Lemma~\ref{HypPol-Lem}. 
\end{itemize}

\begin{rem}
There is another argumentation for cubics: The Fano variety
$F(V_3)$
of the lines on the cubic
$V_3$
is a smooth variety with a very ample canonical 
divisor (\cite[7.8,10.13]{CG72}). 
Furthermore, through a general point there are
exactly 
$6$
lines (\cite{Tju72}).

But as already shown in the beginning of the section, the 
existence of a vector field with zeroes implies an effective linear algebraic
group operation with fixed points on
$V_3$. 
Every group operation on
$V_3$
induces a group operation on
$F(V_3)$, 
and this operation must be trivial because of the very ample 
canonical divisor. Therefore, lines on
$V_3$
are stabilized by the operation, and the intersection point of 
$6$
lines is fixed. The operation on
$V_3$
is trivial, too.
\end{rem}

\subsection{The projective space $\boldsymbol{\mathbb{P}^3}$}

A vector field vanishing on a line is given by
$D_v = x_2\frac{\partial}{\partial x_2} + 
        x_3\frac{\partial}{\partial x_3}$.

\end{document}